\documentclass[11pt]{article}

\usepackage{graphicx,mathtools}
\usepackage{tikz,enumerate, multicol}
\usepackage{float,bbm}
\usepackage{scalerel}
\usepackage{amssymb}
\usepackage{amsfonts}
\usepackage{latexsym}
\usepackage{epsfig}    
\usepackage{amsmath}
\usepackage{amsthm}
\usepackage{graphicx}
\usepackage{subfig}
\usepackage{grffile}

\usepackage[utf8]{inputenc}
\usepackage[english]{babel}

\usetikzlibrary{automata,positioning}
\usetikzlibrary{shapes.geometric,arrows,shapes.misc}

\newtheorem{ittheorem}{Theorem}
\newtheorem{itlemma}{Lemma}
\newtheorem{itdefinition}{Definition}
\newtheorem{itremark}{Remark}

\newenvironment{theorem}{\addtocounter{counter}{1}
\begin{ittheorem}}{\end{ittheorem}}

\newenvironment{lemma}{\addtocounter{counter}{1}
\begin{itlemma}}{\end{itlemma}}

\newenvironment{definition}{\addtocounter{counter}{1}
\begin{itdefinition}}{\end{itdefinition}}

\newenvironment{remark}{\addtocounter{counter}{1}
\begin{itremark}}{\end{itremark}}

\newcommand{\be}[1]{\begin{equation}\label{#1}}
\newcommand{\ee}{\end{equation}}

\newcommand{\bl}[1]{\begin{lemma}\label{#1}}
\newcommand{\el}{\end{lemma}}

\newcommand{\br}[1]{\begin{remark}\label{#1}}
\newcommand{\er}{\end{remark}}

\newcommand{\bt}[1]{\begin{theorem}\label{#1}}
\newcommand{\et}{\end{theorem}}

\newcommand{\bd}[1]{\begin{definition}\label{#1}}
\newcommand{\ed}{\end{definition}}

\newcommand{\bpr}{\begin{proof}}
\newcommand{\epr}{\end{proof}}

\newcommand{\bprl}[1]{\begin{proofof}{\it\ref{#1}}.\,\,}
\newcommand{\eprl}{\end{proofof}}

\newcommand{\bi}{\begin{itemize}}
\newcommand{\ei}{\end{itemize}}

\newcommand{\ben}{\begin{enumerate}}
\newcommand{\een}{\end{enumerate}}

\makeatletter
 \@addtoreset{equation}{section}
 \makeatother

 \makeatletter
 \@addtoreset{counter}{section}
 \makeatother

 \newcounter{counter}[section]

\newcommand{\FF}{\mathcal{F}}
\newcommand{\R}{\mathbb{R}}
\newcommand{\1}{\mathbbm{1}}

\newcommand{\N}{\mathbb{N}}

\newcommand{\PP}{\mathbb{P}}
\newcommand{\EE}{\mathbb{E}}

\newcommand{\eee}{\mathrm{e}}
\newcommand{\PAT}{\mathrm{PAT}}
\newcommand{\SRPAT}{\mathrm{SRPAT}}
 
\begin{document}


\title{Self-Reinforced Preferential Attachment}

\author{\renewcommand{\thefootnote}{\arabic{footnote}}
Yogesh Dahiya
\footnotemark[1]
\\
\renewcommand{\thefootnote}{\arabic{footnote}}
Frank den Hollander
\footnotemark[2]
}

\footnotetext[1]{
Department of Mathematical Sciences, IISER Mohali, Knowledge City, Sector 81, Manauli, PO, Sahibzada Ajit Singh Nagar, Punjab, 140306, India.\\ 
{\sc ph20002@iiser.mohali.ac.in}
}

\footnotetext[2]{
Mathematical Institute, Leiden University, Einsteinweg 55, 2333 CC Leiden, The Netherlands.\\
{\sc denholla@math.leidenuniv.nl}
}

\maketitle

\begin{abstract}
We consider a preferential attachment random graph with self-reinforcement. Each time a new vertex comes in, it attaches itself to an old vertex with a probability that is proportional to the sum of the degrees of that old vertex at all prior times. The resulting growing graph is a random tree whose vertices have degrees that grow polynomially fast in time. We compute the growth exponent, show that it is strictly larger than the growth exponent in the absence of self-reinforcement, and develop insight into how the self-reinforcement affects the growth. Proofs are based on a stochastic approximation scheme.

\medskip\noindent
{\it AMS} 2020 {\it subject classifications.}
05C80, 
60C05, 
60F15. 


\medskip\noindent
{\it Acknowledgment.} YD was supported by an IISER Mohali doctoral fellowship. FdH was supported by the Netherlands Organisation for Scientific Research (NWO) through Gravitation-grant NETWORKS-024.002.003, and by the National Science Foundation (NSF) under Grant No.\ DMS-1928930 while in residence at the Simons Laufer Mathematical Sciences Institute in Berkeley, California, USA during the Spring 2025 semester. The authors thank the International Center for Theoretical Sciences (ICTS) in Bangalore and the Chennai Mathematical Institute (CMI) in Chennai for hospitality in September and December of 2024, and are grateful to Siva Athreya, Vivek Borkar, Remco van der Hofstad and Neeraja Sahasrabudhe for discussions and suggestions.
\end{abstract}


\section{Introduction and main results}
\label{s.intro}


\subsection{Preferential attachment}
\label{ss.PAM}

The \emph{Preferential Attachment Tree} (PAT) is a dynamically growing tree in which vertices arrive one by one and attach themselves to vertices that are already present. Each new vertex connects to a randomly chosen old vertex, such that the probability to connect to a specific vertex is proportional to the current degree of that vertex. The PAT is a random tree sequence $(\PAT_t^\delta)_{t \geq 0}$ whose evolution depends on a parameter $\delta > -1$. At time $t$ there are $t+1$ vertices and $t$ edges in the tree. The following is taken from \cite[Chapter 8]{vdH2017}.

Start at time $t = 0$ with an isolated vertex $v_0$. At time $t = 1$ a vertex $v_1$ comes in and attaches itself to vertex $v _0$ via a single edge. At subsequent times the growth is as follows. For $t \geq 1$, suppose that $v_0,\ldots,v_t$ are the vertices of $\PAT_t^\delta$ with degrees $d_0(t),\ldots,d_t(t)$, respectively. Then, given $\PAT_t^\delta$, the rule for obtaining $\PAT_{t+1}^\delta$ is that a single vertex $v_{t+1}$  attaches itself via a single edge to one of the vertices $v_0,\ldots,v_t$ with probabilities
\[
\PP\big(t+1 \to  i \mid \PAT_t^\delta\big) 
= \dfrac{d_i(t)+\delta}{2t+\delta(t+1)}, \qquad 0 \leq i \leq t,
\]
where $\PP$ denotes the law of $(\PAT_t^\delta)_{t \geq 0}$, $\{t+1 \to i\}$ denotes the event that vertex $v_{t+1}$ attaches itself to vertex $v_i$, and the normalisation uses that $\sum_{i=0}^t d_i(t) = 2t$. 

The role of the parameter $\delta$ is that of a shift: for $\delta = 0$ the attachment probabilities are \emph{linear} in the degrees, for $\delta \neq 0$ they are \emph{affine} in the degrees. The parameter $\delta$ affects the growth of the degrees and the limiting degree distribution. In particular, for every fixed $i \geq 0$, 
\[
d_i(t) \sim \epsilon_i\, t^{1/(\tau-1)}, \qquad t \to \infty, \qquad \PP\text{-a.s.},
\] 
with $\tau = 3+\delta$ a characteristic \emph{exponent} and $\epsilon_i \in (0,\infty)$ a random variable \cite[Section 8.3, Theorem 8.14 and Lemma 8.17]{vdH2017}. Moreover, for every fixed $k \geq 1$,
\[
\lim_{t\to\infty} \tfrac{1}{t+1} \sum_{i=0}^t 1_{\{d_i(t)=k\}} = p_k, \qquad \PP\text{-a.s.},
\] 
with $p_k$ given by a quotient of Gamma functions, and $p_k \sim c_\delta k^{-\tau}$ as $k\to\infty$ with prefactor $c_\delta \in (0,\infty)$ \cite[Theorem 8.3]{vdH2017}. 

Other versions of preferential attachment random graphs have been considered as well, allowing for self-loops and multiple edges. Their scaling behaviour is similar \cite[Sections 8.2--8.3]{vdH2017}. 


\subsection{Self-reinforced preferential attachment}
\label{ss.SRPAM}

In the present paper we consider a version of the PAT with \emph{self-reinforcement} and shift $\delta=0$, which we refer to as SRPAT. The way in which the self-reinforcement is implemented is similar to what is done for \emph{self-reinforced random walk} on a graph, where the probability to cross an edge is proportional to the total number of crossings of that edge at all prior times \cite{MR2006}. 

Start at time $t=0$ with an isolated vertex $v_0$. At time $t=1$ a vertex $v_1$ comes in and attaches itself to vertex $v_0$ via a single edge. At subsequent times the growth is as follows. At time $t \geq 1$ a vertex $v_{t+1}$ enters and attaches itself via a single edge to one of the vertices already present in the graph with a probability that is proportional to the sum of the degrees of that old vertex at all prior times, i.e.,
\[
\PP\big(t+1 \to i \mid (\SRPAT_t)_{0 \leq s \leq t}\big) 
= \frac{\theta_t(i)}{\sum_{j=0}^t \theta_t(j)}, \qquad 0 \leq i \leq t,
\]
with
\begin{equation}
\label{eq:weightchoice}
\theta_t(i) = \sum_{s=1}^t d_s(i), \qquad 0 \leq i \leq t,
\end{equation}
where $d_s(i)$ is the degree of vertex $v_i$ at time $s$, and $\PP$ denotes the law of $(\SRPAT_t)_{t \geq 0}$. More formally, let $\{\FF_t\}_{t \geq 1}$ be the filtration defined by 
\[
\FF_t = \sigma\big(d_s(i)\colon\,1 \leq s \leq t,\, 0 \leq i \leq t\big).
\] 
Then the conditional probability above refers to conditioning on $\FF_t$. As long as only the degrees of the SRPAT are monitored, this is the proper filtration to work with. Note that $d_1(0)=d_1(1)=1$, while, for $i \geq 2$, $d_s(i) = 0$, $1 \leq s < i$, and $d_i(i) = 1$. 

Clearly, for every $t \geq 1$,
\begin{equation}
\label{eq:rec}
\begin{array}{lll} 
d_{t+1}(i) &=& d_t(i) + \1_{\{t+1 \sim i\}},\\[0.2cm]
\theta_{t+1}(i) &=& \theta_t(i) + d_t(i) + \1_{\{t+1 \sim i\}}, 
\end{array}
\qquad 0 \leq i \leq t.
\end{equation}
Since there are $s$ edges in the graph at time $s$, we have $\sum_{i=0}^t \theta_t(i) = 2 \sum_{s=1}^t s = t(t+1)$, $t \geq 1$. Hence, for every $t \geq 1$,
\begin{equation}
\label{eq:law} 
\PP(t+1 \to i \mid \FF_t) =  \frac{\theta_t(i)}{t(t+1)}, \qquad 0 \leq i \leq t.
\end{equation}
It follows from \eqref{eq:rec}--\eqref{eq:law} that, for every $t \geq 1$,
\begin{equation}
\label{eq:meanrec}
\begin{array}{lll}
\EE[d_{t+1}(i) \mid \FF_t] &=& d_{t}(i) + \frac{\theta_t(i)}{t(t+1)},\\[0.2cm]
\EE[\theta_{t+1}(i) \mid \FF_t] &=& d_t(i) + \theta_t(i)\left(1 + \frac{1}{t(t+1)}\right),
\end{array}
 \qquad 0 \leq i \leq t.
\end{equation}


\subsection{Theorems}
\label{ss.theorems}

Our first two theorems provide a comparison of the weights and the degrees. Since $v_0$ and $v_1$ are stochastically symmetric we can focus on $i \geq 1$. For $t\geq i$, define
\begin{equation}
\label{eq:alphabeta}
\alpha_t(i) = \frac{t\,d_t(i)}{\theta_t(i)}, \qquad \beta_t(i) = \frac{t\,\EE[d_t(i)]}{\EE[\theta_t(i)]},
\qquad 0 \leq i \leq t,
\end{equation}
and note that $\alpha_i(i) = \beta_i(i) = i$.

\begin{theorem}
\label{fixed_point_convergence} 
For every $i\geq 1$, $\lim_{t\to\infty} \alpha_t(i) = \phi$ $\PP$-a.s.\ with $\phi = \tfrac12(1+\sqrt{5}) \approx 1.618$. 
\end{theorem}

\begin{theorem} 
\label{mode_convergence}
For every $i \geq 1$, $\lim_{t \to \infty} \beta_t(i) = \phi$. Moreover, there exists a finite $T(i) \geq i$ such that $t \mapsto \beta_t(i)$ is strictly decreasing for $i \leq t < T(i)$ and strictly increasing for $t \geq T(i)$.
\end{theorem}

\noindent
Our main theorem identifies the growth of the degrees.

\begin{theorem}
\label{degree_growth} 
For every $i\geq 1$, $\lim_{t\to\infty} t^{-1/\phi} d_t(i) = \epsilon_i$ $\PP$-a.s.\ with $\EE[\epsilon_i] \in (0,\infty)$. Moreover, $\lim_{i \to \infty} i^{1/\phi} \EE[\epsilon_i] = 1$. 
\end{theorem}

Theorems~\ref{fixed_point_convergence} and \ref{degree_growth} are proved in Section~\ref{s.proofs} with the help of a stochastic approximation scheme, Theorem~\ref{mode_convergence} with the help of an iterated mapping scheme. The proofs also provide rate of convergence estimates. Appendix \ref{appA} recalls what a stochastic approximation scheme is and under what conditions it can be employed.


\subsection{Discussion}

Our paper is the first to consider preferential attachment with self-reinforcement. The self-reinforcement introduces \emph{memory} into the growth of the random tree, which needs to be harnessed.  

\vspace{-0.3cm}
\begin{figure}[htbp]
\begin{center}
\includegraphics[scale=0.35]{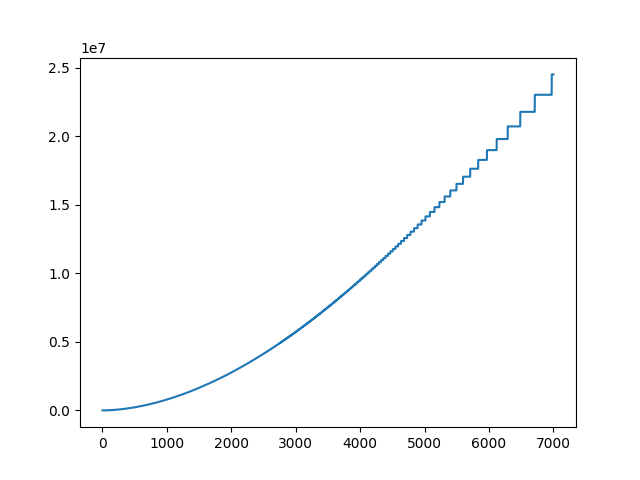}
\end{center}
\vspace{-0.5cm}
\caption{\small Plot of $i \mapsto T(i)$ with horizontal scale $10^3$ and vertical scale $10^7$.}
\label{fig:vertex}
\end{figure}  

\begin{itemize}
\item
Theorem~\ref{fixed_point_convergence}, which says that $\theta_t(i) \sim \frac{1}{\phi} t d_t(i)$ as $t\to\infty$, shows that the \emph{quenched version} of the SRPAT behaves as the quenched version of the PAT when in the latter the probability for the new vertex at time $t$ to attach to an old vertex is asymptotically proportional to the degree of the old vertex divided by $\phi\,t$.
\item
Theorem~\ref{mode_convergence}, which says that $\EE[\theta_t(i)] \sim \frac{1}{\phi} t\, \EE[d_t(i)]$ as $t\to\infty$, shows that the \emph{annealed version} of the SRPAT is close to the quenched version. It also gives information on the mode of convergence in the annealed version. Figure~\ref{fig:vertex} shows a plot of the crossover time $T(i)$ as a function of $i$.
\item
Theorem~\ref{degree_growth} shows that the degrees of the vertices grow like $t^{1/\phi}$, which is faster than $t^{1/2}$ for the PAT with shift $\delta = 0$ because $\phi<2$. The exponent $1/\phi$ for the SRPAT corresponds to the exponent $1/(2+\delta)$ for the PAT with shift $\delta = \phi - 2 \approx - 0.382$. Since $\EE[\epsilon_i] \in (0,\infty)$, we know that $\epsilon_i<\infty$ $\PP$-a.s. We expect that $\epsilon_i>0$ $\PP$-a.s.\ (simulations support this expectation), but are unable to prove this. The degree sequence at vertex $i$ does not appear to be tail trivial.
\item
Here is a consistency check for the growth exponent in Theorem~\ref{degree_growth}. If $d_t(i) \sim \epsilon_i\,t^{1/\phi}$, then 
\[
\theta_t(i) = \sum_{s=1}^t d_s(i) \sim \epsilon_i \sum_{s=1}^t s^{1/\phi} \sim \epsilon_i\,\frac{\phi}{\phi+1}\, t^{(\phi+1)/\phi}.
\] 
Combined with $td_t(i) \sim \epsilon_i\,t^{(\phi+1)/\phi}$, this gives that $\alpha_t(i)$ converges to $(\phi+1)/\phi$. The latter matches with Theorem~\ref{fixed_point_convergence}, because $\phi$ is the unique positive solution to the equation $(\phi+1)/\phi = \phi$, i.e., $\phi^2-\phi-1=0$.
\item  
The classical preferential attachment tree (PAT) described in Section~\ref{ss.PAM} assumes that the probability of attachment is a linear function of the degree. This linear form is known to produce scale-free networks with power-law degree distributions. Extensions of this model have been studied where instead the attachment probability is proportional to $d_t(i)^\alpha$ for $\alpha > 0$ \cite{dereich2009random,sethuraman2016growth}. In the \emph{sublinear} case $\alpha < 1$, the resulting graphs fail to exhibit power-law behavior, while in the \emph{superlinear} case $\alpha > 1$, a condensation phenomenon occurs in which a few nodes acquire a macroscopic fraction of all edges. Only the linear case $\alpha = 1$ yields networks with true scale-free structure.
\item
Our paper demonstrates that the class of preferential attachment models producing power-law behaviour is larger. In our model the attachment probability scales like $\mathbb{P}(t+1 \to i \mid \mathcal{F}_t) \sim \frac{1}{\phi} \frac{d_t(i)}{t}$ a.s.\ as $t \to \infty$. We hypothesise that a power-law degree distribution arises whenever the attachment probability scales like $\mathbb{P}(t+1 \to i \mid \mathcal{F}_t) \sim c\frac{d_t(i)}{t}$ a.s.\ as $t\to\infty$ for some constant $c > 0$, with the power-law exponent determined by the value of $c$. In our model, $c = 1/\phi$, and  the resulting degree distribution is expected to follow a power law with exponent $1 + \phi$ (see Section~\ref{ss.PAM}). 
\end{itemize}


\section{Proofs}
\label{s.proofs}

In Sections \ref{ss.fixedpointconv}--\ref{ss.rate} we give the proof of Theorems~\ref{fixed_point_convergence}--\ref{degree_growth}, respectively. Along the way we need a technical lemma, which we state and prove in Section~\ref{ss.lem}.


\subsection{A priori growth of degrees}
\label{ss.lem}

For $i \geq 1$, let
\begin{equation}
\label{eq:eventdef}
\mathcal{E}(i) = \bigcap_{k\in\N} \left\{\liminf_{t\to\infty} d_t(i)/(\log t)^k \geq 1/k!\right\}.
\end{equation}

\begin{lemma}
\label{lem:growth} 
For every $i \geq 1$, $\PP(\mathcal{E}(i)) = 1$.
\end{lemma}

\begin{proof}
The proof uses the following elementary fact. Let $\{X_u\}_{u \geq 0}$ be a sequence of independent Bernoulli random variables with $\PP(X_u=1) = p_u$ satisfying $\sum_{u \geq 0} p_u = \infty$. Then 
\[
\lim_{v\to\infty} \frac{\sum_{u=0}^v \1_{\{X_u=1\}}}{\sum_{u=0}^v p_u} = 1 \qquad \PP\text{-a.s.}
\]  
Indeed, the weak law of large numbers is immediate from Chebyshev's inequality, while the strong law of large numbers follows from a subsequence argument for exceedences and Borel-Cantelli.

Fix $i \geq 1$. Consider vertex $v_i$, which arrives at time $i$. Since each vertex is reinforced by its incident edges, we have that
\[
\PP\big(s+1 \to i \mid (\SRPAT_u)_{0 \leq u \leq s}\big) \geq \frac{s-i +1}{s(s+1)}, \qquad s \geq i,
\] 
where the bound is \emph{uniform} in the realisation of $(\SRPAT_u)_{0 \leq u \leq s}$. Consequently, there exists an \emph{ordered coupling} of the sequence $\{d_t(i)\}_{t \geq i}$ and the sequence $\{\sum_{0 \leq u \leq v} \1_{\{X_u=1\}}\}_{v \geq 0}$ for the choice $p_u = (u+1)/(u+i)(u+i+1)$ (after time is shifted by $i$). Since for this choice $\sum_{u=0}^v p_u \sim \log v$ as $v\to\infty$, it follows that $\liminf_{t\to\infty} d_t(i)/\log t \geq 1$ $\PP$-a.s. 

Next, we boost up the above estimate. Since
\[
\PP\big(s+1 \to i \mid (\SRPAT_u)_{0 \leq u \leq s}\big) = \frac{\sum_{u=i}^s d_u(i)}{s(s+1)}, \qquad s \geq i,
\] 
and $\liminf_{s\to\infty} \sum_{u=i}^s d_u(i)/(s \log s) \geq 1$ $\PP$-a.s.\ by the lower bound just established, we can repeat the coupling, this time with $p_u \sim (u \log u)/(u+i)(u+i+1)$ as $u\to\infty$. Since for this choice $\sum_{u=0}^v p_u \sim \tfrac12 (\log v)^2$ as $v\to\infty$, it follows that $\liminf_{t\to\infty} d_t(i)/(\log t)^2 \geq \tfrac12$ $\PP$-a.s. 

The same argument can be repeated to get that $\liminf_{t\to\infty} d_t(i)/(\log t)^k \geq 1/k!$ $\PP$-a.s.\ for any $k\in\N$. Since the latter events are monotone in $k$, the claim follows.
\end{proof}


\subsection{Proof of Theorem~\ref{fixed_point_convergence}}
\label{ss.fixedpointconv}

\begin{proof}
We use the \emph{stochastic approximation scheme} outlined in \cite[Chapter 2]{B2023} (see Appendix \ref{appA}). Fix $i \geq 1$ and suppress it from the notation. Put $\alpha^*_t = 1/\alpha_t$ and note that $\tfrac{1}{t} \leq \alpha^*_t \leq 1$ by \eqref{eq:weightchoice} and \eqref{eq:alphabeta}, because all the degrees are non-decreasing in time. Recall that $\alpha^*_i = 1/i$. Define the \emph{martingale difference}
\begin{equation}
\label{eq:martdiff}
M^*_{t+1} = t\alpha^*_t\left(\frac{\1_{\{t+1 \to i\}}}{d_t+\1_{\{t+1 \to i\}}} 
- \EE\left[ \frac{\1_{\{t+1 \to i\}}}{d_t+\1_{\{t+1 \to i\}}} ~\Big|~ \FF_t\right]\right),
\qquad t \geq i.
\end{equation}
Compute, using \eqref{eq:rec} and \eqref{eq:law},
\begin{equation}
\label{eq:alpharec}
\begin{aligned}
\alpha^*_{t+1}  
&= \frac{\theta_{t+1}}{(t+1)d_{t+1}} = \frac{\theta_t + d_t + \1_{\{t+1 \to i\}}}{(t+1)(d_t + \1_{\{t+1 \to i\}})}
= \frac{1}{t+1}\left(1 + t\alpha^*_t\,\frac{d_t}{d_t + \1_{\{t+1 \to i\}}} \right)\\
&=\alpha^*_t + \frac{1}{t+1}\left(1-\alpha^*_t-M^*_{t+1}-{\alpha^*_t}^2\frac{td_t}{(t+1)(d_t+1)}\right)\\
&=\alpha^*_t + \frac{1}{t+1}\left(1-\alpha^*_t- {\alpha^*_t}^2-M^*_{t+1} \right) 
- \frac{{\alpha^*_t}^2}{(t+1)} \left(\frac{t{d_t}}{(t+1)(d_t+1)}-1\right)
\end{aligned}
\end{equation}
and
\[
\begin{aligned}
\EE\big[(M^*_{t+1})^2  \mid \FF_t\big]
&= t^2 {\alpha^*_t}^2\,\EE\left[\left(\frac{\1_{\{t+1 \to i\}}}{d_t
+\1_{\{t+1 \to i\}}}-\frac{\theta_t}{t(t+1)(d_t+1)}\right)^2 ~\Big|~\FF_t \right]\\
&= \alpha^{*2}_t \frac{\theta_t}{(t+1)^2(d_t+1)^2} \big(t(t+1)-\theta_t\big),
\end{aligned}
\]
where we note that $0 \leq \theta_t \leq t(t+1)$. Put 
\begin{equation}
\label{eq:zetadef}
\zeta_t = \sum_{s =i}^{t-1} \frac{M^*_{s+1} }{s+1}, \qquad t \geq i,
\end{equation}
and compute
\[
\begin{aligned}
&\sum_{t \geq i} \EE\big[{|\zeta_{t+1}-\zeta_t|}^2  \mid \FF_t\big]
= \sum_{t \geq i} \frac{1}{(t+1)^2}\,\EE\big[(M^*_{t+1})^2  \mid \FF_t\big]\\
&= \sum_{t \geq i} \alpha^{*2}_t\theta_t\,\frac{t}{(t+1)^3(d_t+1)^2}
= \sum_{t \geq i} \alpha^{*3}_t\,\frac{t^2d_t}{(t+1)^3(d_t+1)^2} \leq \sum_{t \geq i} \frac{1}{t d_t},
\end{aligned}
\]
where we use that $\alpha^*_t \in (0,1]$ for all $t \geq i$. 

On the event $\mathcal{E}(i)$, which has probability $1$ by Lemma \ref{lem:growth}, we have $\sum_{t \geq i} (1/t d_t)<\infty$, and so the sequence $\{\zeta_t\}_{t \geq i}$ is a \emph{zero-mean square-integrable martingale}. Now apply Lemma~\ref{lem:sas} in Appendix~\ref{appA} with
\[
x_t = \alpha^*_t, \quad a_t = \frac{1}{t+1}, \quad h(x) = 1-x-x^2, \quad M_{t+1} = - M^*_{t+1},
\quad \epsilon_{t+1} = \alpha^{*2}_t \left(1-\frac{t}{t+1} \frac{d_t}{d_t+1}\right)
\]
to the recursion in \eqref{eq:alpharec} and note that, on the event $\mathcal{E}(i)$, $\lim_{t\to\infty} td_t/(t+1)(d_t+1)=1$ because $d_t \to\infty$. The conditions needed for Lemma~\ref{lem:sas} are all fulfilled for this choice of variables, and so it follows that $\alpha^*_t$ converges $\PP$-a.s.\ as $t\to\infty$ to the unique positive zero of the quadratic equation $1-x-x^2=0$, which equals $1/\phi$ and is globally attracting in $[0,1]$.
\end{proof}


\subsection{Rate of convergence estimate}
 
There is a large body of literature that establishes rates of convergence for the stochastic approximation scheme given in Appendix~\ref{appA} (see e.g.\ \cite{R1982}, \cite{K1998}, \cite{P1999}, \cite{KV2024}). However, in our setting the variance of the martingale difference sequence $\{M^*_t\}_{t \geq i}$ in \eqref{eq:martdiff} diverges, and so such a refinement is non-standard. Nonetheless, by using Lemma \ref{lem:sasrate} we are able to get the following. 

\begin{lemma}
\label{lem:alphastarconv}
$\lim_{t\to\infty} t^\lambda|\alpha^*_t-(1/\phi)| = 0$ $\PP$-a.s.\ for all $\lambda \in (0,\infty)$.
\end{lemma}

\begin{proof}
The proof comes in four steps.

\medskip\noindent
{\bf 1.} We begin with a fluctuation bound for the martingale defined in \eqref{eq:zetadef}. 

\begin{lemma}
\label{lem:martest}
$\lim_{t \to \infty} (\log t)^\ell \sup_{t \le s \le 2t} |\zeta_s - \zeta_t| = 0$ $\PP$-a.s.\ for all $\ell\in\N$.
\end{lemma}

\begin{proof}
Fix $\ell \in \N$. For $n\in\N$, put $t_n = 2^n$ and $X_n = \sup_{t_n \leq s \le t_{n+1}} |\zeta_s - \zeta_{t_n}|$. Since $\{\zeta_t\}_{t \geq i}$ is a square-integrable martingale, Doob's maximal inequality gives
\[
\EE[X_n^2] = \EE\left[\sup_{t_n \leq s \leq 2t_n} |\zeta_s - \zeta_{t_n}|^2\right]
\leq 4\,\EE\big[|\zeta_{2t_n} - \zeta_{t_n}|^2\big].
\]
We saw in Section \ref{ss.lem} that $d_s(i) \geq \frac{1}{k!} (\log s)^k$ $\PP$-a.s.\ for all $k\in \N$ and all $s$ large enough (recall Lemma~\ref{lem:growth}), and in Section \ref{ss.fixedpointconv} that $\EE[(M^*_{s+1})^2 \mid \mathcal{F}_s] \leq s/d_s(i)$ for all $s \geq i$. Since $\zeta_{s+1} - \zeta_s = M^*_{s+1}/(s+1)$, we therefore get, for all $s$ large enough,
\[
\EE[|\zeta_{s+1} - \zeta_s|^2] \leq \frac{C_k}{s (\log s)^k}
\]
for any $C_k \in (k!,\infty)$, and hence
\[
\EE[|\zeta_{2t_n} - \zeta_{t_n}|^2] \leq \sum_{t_n \leq s < 2t_n} \frac{C_k}{s (\log s)^k}.
\]
Since $\log s \geq n \log 2$ for all $s \geq t_n$, we get, for all $n$ large enough, $\EE[X_n^2] \leq \frac{C_k'}{n^k}$ for any $C_k' \in (4k!/(\log 2)^{k-1},\infty)$. Next, put $\varepsilon_n = n^{-\ell}$. Then, by Markov's inequality,
\[
\PP(X_n > \varepsilon_n) \leq \frac{\EE[X_n^2]}{\varepsilon_n^2}
\leq \frac{C_k'}{n^{k-2\ell}}.
\]
Define the event $A_n = \{X_n > \varepsilon_n\}$. Pick $k > 2\ell+1$. Then $\sum_{n\in\N} \PP(A_n) < \infty$ and so, by the Borel-Cantelli lemma, $\PP(A_n \text{ infinitely often}) = 0$, i.e., $X_n \le \frac{1}{n^\ell}$ for $n$ large enough $\PP$-a.s. Finally, if $t \in [t_n, t_{n+1}]$, then $[t, 2t] \subset [t_n,t_{n+2}]$ and so, by the triangle inequality,
\[
\begin{aligned}
&\sup_{t \leq s \le 2t} |\zeta_s - \zeta_t| 
\leq \sup_{t \leq s \le 2t} |\zeta_s - \zeta_{t_n}| + |\zeta_t-\zeta_{t_n}| 
\leq 2\sup_{t_n \leq s \leq t_{n+2}} |\zeta_s - \zeta_{t_n}|\\
&\leq 2\sup_{t_n \leq s \le t_{n+1}} |\zeta_s - \zeta_{t_n}| 
+ 2\sup_{t_{n+1} \leq s \le t_{n+2}} |\zeta_s - \zeta_{t_{n+1}}|
= 2X_n + 2X_{n+1} \leq \frac{2}{n^\ell} + \frac{2}{(n+1)^\ell},
\end{aligned}
\]
where the last inequality holds for $n$ large enough $\PP$-a.s. Since $n+1 = \log t_{n+1}/\log 2 \geq \log t/\log 2$, and $\ell \in \N$ is arbitrary, this yields the claim.
\end{proof}

\medskip\noindent
{\bf 2.}
We apply Lemma \ref{lem:sasrate} for the choice $T = T_{2t}-T_t$. In our particular setting we have $h(0)=1$, $L=3$, $C_0=1$. Moreover, as $t\to\infty$,
\begin{equation}
\label{eq:estcoef}
\sum_{t \leq s \leq 2t} a^2_s = \sum_{t \leq s \leq 2t} (s+1)^{-2} = \tfrac12 t^{-1} + O(t^{-2}), 
\quad \sup_{t \leq s \leq 2t} a_s = a_t = (t+1)^{-1},\\ 
\sum_{t \leq s \leq 2t} a_s |\epsilon_s|  = e_\PP(t),
\end{equation}
where we introduce $e_\PP(t)$ as short-hand notation for any function that $\PP$-a.s.\ decays faster than any power of $1/\log t$ (recall Lemma~\ref{ss.lem}). The bound in \eqref{comparison} says that, for all $t \geq 0$, 
\begin{equation}
\label{est1}
\sup_{u \in [T_t,T_{2t}]} |\bar{\alpha}(u) - y^{T_t}(u)| 
\leq K_t\,\eee^{3(T_{2t}-T_t)} + C_t\,(t+1)^{-1} + O\big(e_\PP(t)\big) \qquad \PP\text{-a.s.},
\end{equation}
where $(\bar{\alpha}(u))_{u \in [0,\infty)}$ is the piecewise linear interpolation of $\{\alpha^*_t\}_{t \geq 0}$ over the successive time intervals $u \in [T_t,T_{t+1}]$, $t \geq 0$, and $(y^{T_t}(u))_{u \in [T_t,\infty)}$ is the solution to the ODE $y'(u) = h(y(u))$, $u \geq 0$, with initial value $y^{T_t}(T_t) = \alpha^*_{T_t}$, and
\[
\begin{aligned}
K_t &= 3C_t\left(\tfrac12 t^{-1}+O(t^{-2})\right) + O\big(e_\PP(t)\big), \qquad t \to \infty,\\ 
C_t &= 1+3(1+T_{2t}-T_t)\,\eee^{3(T_{2t}-T_t)}.
\end{aligned}
\] 
where the error term in the first line comes from Lemma \ref{lem:martest} and \eqref{eq:estcoef}.

\medskip\noindent
{\bf 3.}
Since $(y^{T_t}(u)-\psi)' = h(y^{T_t}(u)) - h(\psi) = -(1+2\psi)(y^{T_t}(u)-\psi) - (y^{T_t}(u)-\psi)^2$, we have
\[
|y^{T_t}(u)-\psi|' = -(1+2\psi)|y^{T_t}(u)-\psi| - |y^{T_t}(u)-\psi|^2\, \text{sign}(y^{T_t}(u)-\psi).
\]
Since $|y^{T_t}(u)-\psi| \leq 1$, it follows that
\[
-2\psi\, |y^{T_t}(u)-\psi| \geq |y^{T_t}(u)-\psi|'
\] 
and hence
\begin{equation}
\label{est2}
|y^{T_t}(u)-\psi| \leq |\alpha^*_{T_t}-\psi|\,\eee^{-2\psi\,(u-T_t)}, \qquad u \in [T_t,T_{2t}].
\end{equation}

\medskip\noindent
{\bf 4.}
Since $T_t = \sum_{s=0}^{t-1} (s+1)^{-1} = \log t + \gamma + O(t^{-1})$, with $\gamma$ Euler's constant, we have $T_{2t}-T_t = \log 2 + O(t^{-1})$. From \eqref{est1}--\eqref{est2} we get that
\begin{equation}
\label{est3}
\begin{aligned}
|\alpha^*_{T_{2t}}-\psi| &\leq |\alpha^*_{T_{2t}}- y^{T_t}(T_{2t})| + |y^{T_t}(T_{2t})-\psi|\\
&\leq \big[At^{-1} + O(t^{-2}) + O\big(e_\PP(t)\big)\big] + |\alpha^*_{T_t}-\psi|\,\eee^{-2\psi\,[\log 2+O(t^{-1})]}
\end{aligned}
\end{equation}
with $A = 325+312\log 2$. Abbreviate $\delta_t = |\alpha^*_{T_t}-\psi|$. Then the estimate in \eqref{est3} reads
\[
\delta_{2t} \leq O\big(e_\PP(t)\big) + \delta_t B_t \quad \text{ with } 
\quad B_t = \eee^{-2\psi[\log 2+O(t^{-1})]}. 
\]  
By Theorem \ref{fixed_point_convergence}, for every $\epsilon>0$ $\PP$-a.s.\ there exists a random $\mathcal{T}_\epsilon<\infty$ such that $\delta_t \leq \epsilon$ for all $t \geq \mathcal{T}_\epsilon$. Next, fix $t \geq \mathcal{T}_\epsilon$ such that $B_t \leq B = \eee^{-\psi \log 2} = 2^{-\psi}$. Iterate the last inequality to get, for $t \geq \mathcal{T}_\epsilon$, 
\[
\delta_{2^{n+1}t} \leq \sum_{m=0}^n B^m \, O\big(e_\PP(2^{n-m}t)\big) + B^{n+1}\,\epsilon, \qquad n \in \N_0.
\] 
Since $B < 1$, the terms with small $m$ dominate and we get, for $t \geq \mathcal{T}_\epsilon$,  
\[
|\alpha^*_{T_{2^{n+1}t}}-\psi| = O\big(e_\PP(2^nt)\big), \qquad n \to\infty.
\]
Since $T_{2^{n+1}t} \sim \log (2^{n+1}t)$ as $n\to\infty$, this is the same as saying that 
\[
|\alpha^*_{T_{2^{n+1}t}}-\psi| = O\big(e_\PP(\eee^{T_{2^nt}})\big), \qquad n \to \infty.
\]
The latter settles the claim along the time subsequence $\{T_{2^nt}\}_{n\in\N}$ for fixed $t\geq \mathcal{T}_\epsilon$, and hence for all times by interpolation.
\end{proof} 


\subsection{Proof of Theorem~\ref{mode_convergence}}
\label{ss.modeconv}

\begin{proof}
Note that $1 \leq \beta_t \leq t$ by \eqref{eq:weightchoice} and \eqref{eq:alphabeta}, because all the degrees are non-decreasing in time. Recall that $\beta_i=i$. As before, fix $i \geq 1$ and suppress it from the notation. By \eqref{eq:meanrec}--\eqref{eq:alphabeta}, we have 
\[
\beta_{t+1} = F_t(\beta_t), \qquad t \geq i,
\] 
where 
\[
F_t(x)=\frac{x(t+1)+1}{x+t+\frac{1}{t+1}}.
\] 
The fixed point equation $F_t(x)=x$ reads $x^2-\frac{t}{t+1}\,x-1=0$, which gives 
\[
x_t=\frac{1}{2}\left(\frac{t}{t+1} +\sqrt{{\left(\frac{t}{t+1}\right)}^2+4}\,\right).
\] 
Note that $t \mapsto x_t$ is strictly increasing with $\lim_{t\to\infty} x_t = \phi$.

Compute 
\[
\begin{aligned}
F'_t(x) &= \frac{t(t+1)^3}{(x(t+1)+t(t+1)+1)^2} > 0,\\
F''_t(x) &= \frac{-2t(t+1)^4}{(x(t+1)+t(t+1)+1)^2} < 0, 
\end{aligned}
\]
i.e., for every $t \geq i$, $x \mapsto F_t(x)$ is strictly increasing and strictly concave. Consequently, $x_t$ is an \emph{attracting fixed point} of $F_t$, in particular, $F_t(x)<x$ for all $x>x_t$ and $F_t(x)>x$ for all  $x<x_t$ (see Figure~\ref{fig:map}). It follows that either $t \mapsto \beta_t$ is decreasing (Case 1) or there exists a finite $T(i)$ such that it is decreasing for $t<T(i)$ and increasing for $t\geq T(i)$ (Case 2). In both cases the sequence $\{\beta_t\}_{t \geq 1}$ has a limit, say $L \in (0,\infty)$. We want to show that $L=\phi$, i.e., $\lim_{t\to\infty} |\beta_t-x_t|=0$. Since $\lim_{t \to \infty} F'_t(x) = 1$ for all $x \in (0,\infty)$, $F_t$ converges to the identity map as $t \to \infty$ (on compacts), and so the latter target is not obvious.

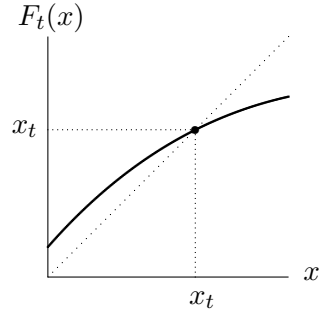
\begin{figure}[htbp]
\vspace{0.5cm}
\begin{center}
\setlength{\unitlength}{0.4cm}
\begin{picture}(8,8)(0,0)
\put(0,0){\line(8,0){8}}
\put(0,0){\line(0,8){8}}
{\thicklines
\qbezier(0,1)(3.5,5)(8,6)
}
\qbezier[50](0,0)(4,4)(8,8)
\qbezier[30](0,4.9)(2.45,4.9)(4.9,4.9)
\qbezier[30](4.9,0)(4.9,2.45)(4.9,4.9)
\put(4.9,4.9){\circle*{.25}}
\put(8.5,-0.1){$x$}
\put(-1,8.4){$F_t(x)$}
\put(-1.2,4.8){$x_t$}
\put(4.7,-0.9){$x_t$}
\end{picture}
\end{center}
\vspace{0.2cm}
\caption{\small The drawn curve is the map $x \mapsto F_t(x)$, the dotted line is the identity map.}
\label{fig:map}
\end{figure}

\medskip\noindent 
\underline{Case 1}: $\beta_t \downarrow L$. Note that, as long as $\beta_t > x_t$, we have $\beta_t - \beta_{t+1} \geq (1 - F'_t(x_t)) (\beta_t - x_t)$ by the concavity of $F_t$. Also note that $1-F'_t(x) = \frac{2x-1}{t} + O(\frac{1}{t^2})$ as $t\to\infty$. Defining $\Delta_t = \beta_t - x_t$, $\eta_t = 1-F'_t(x_t)$ and $\delta_t= x_{t+1}-x_t$ for $t \geq i$, we have
\[
0 \leq \Delta_{t+1} = -\delta_t + \Delta_t + (\beta_{t+1} - \beta_t) 
\leq -\delta_t + (1-\eta_t)\Delta_t \leq  (1-\eta_t)\Delta_t.
\]
Hence $0 \leq \Delta_{t+1} \leq \Delta_i\prod_{s=i}^t(1-\eta_s)$, with the upper bound tending to zero as $t \to \infty$ because $\eta_s  \sim \frac{2\phi-1}{s}$ as $s \to \infty$ and $2\phi-1=\sqrt{5}>0$.

\medskip\noindent 
\underline{Case 2}: $\beta_t \uparrow L$ for $t \geq T(i)$. Recall that $t \mapsto \beta_t$ lies to the right of $1$. Since, by assumption, $\beta_t < x_t$ for $t \geq T(i)$, we have $\beta_{t+1}-\beta_t \geq (1 - F'_t(\beta_t))(x_t -\beta_t)$ for $t \geq T(i)$ by the concavity of $F_t$. Defining $\tilde{\Delta}_t =x_t-\beta_t$ and $\tilde{\eta}_t = F'_t(\beta_t)$, we have 
\[
0 \leq \tilde{\Delta}_{t+1} = \delta_t + \tilde{\Delta}_t - (\beta_{t+1}-\beta_t)
\leq \delta_t + \tilde{\Delta}_t - (1 - \tilde{\eta}_t) \tilde{\Delta}_t
= \delta_t + \tilde{\eta}_t \tilde{\Delta}_t.
\]
Hence 
\[
0 \leq \tilde{\Delta}_t \leq \sum_{s=T(i)}^{t-1} \delta_s \prod_{u=s+1}^{t-1} \tilde{\eta}_u 
+ \tilde{\Delta}_{T(i)} \prod_{s=T(i)}^{t-1} \tilde{\eta}_s.
\]
Since $\tilde{\eta}_s < F'_s(1)$ for all $s$, this gives, for $t \geq T(i)$,
\[
0 \leq \tilde{\Delta}_t \leq \sum_{s=T(i)}^{t-1} \delta_s \prod_{u=s+1}^{t-1} F'_u(1) 
+ \tilde{\Delta}_{T(i)} \prod_{s=T(i)}^{t-1} F'_s(1).
\]
Since $1-F'_s(1) = \frac{1}{s} + O(\frac{1}{s^2})$ as $s \to \infty$, it follows that $\lim_{t\to\infty} \prod_{u=s+1}^{t-1} F'_u(1) = 0$ for all $s \geq T(i)$. Since $\lim_{t\to\infty} \sum_{s=T(i)}^{t-1} \delta_s = \phi-x_{T(i)}$, we get that $\lim_{t\to\infty} \tilde{\Delta}_t = 0$.

\medskip
To complete the proof we show that the scenario in Case 1 never occurs, i.e., eventually $\beta_t$ goes below $x_t$, and stays below $x_t$ while moving up together with $x_t$ to $\phi$. The proof is by contradiction. Suppose that $\beta_t \downarrow \phi$. Then, as shown above,   
\[
\Delta_{t+1} \leq -\delta_t + (1-\eta_t)\Delta_t \leq -\delta_t + \Delta_i \prod_{s=i}^t (1-\eta_s).
\]
It follows from the formula for $x_t$ that $\delta_t \sim \frac{\phi}{2\phi-1}\frac{1}{t^2}$ as $t\to\infty$. On the other hand,  $\prod_{s=i}^t(1-\eta_s) \leq \exp[-\sum_{s=i}^t \eta_s] = \exp[-(2\phi-1) \log t + O(1)] = \Theta(1)\, t^{-(2\phi-1)}$ as $t\to\infty$. Since $2\phi-1=\sqrt{5}>2$, the second term in the right-hand side tends to zero faster than the first term in the right-hand side. Hence $\Delta_{t+1}$ eventually becomes negative, which contradicts the assumption.
\end{proof} 

The following rate of convergence estimate is the analogue of Lemma \ref{lem:alphastarconv}.

\begin{lemma}
\label{lem:betaconv}
$\limsup_{t\to\infty} \frac{t}{\log t} |\beta_t-\phi| < \infty$.
\end{lemma}

\begin{proof}
For $t \geq T(i)$, we have $0 \leq \phi-\beta_t = \phi - x_t + \tilde{\Delta}_t$. We have $\phi - x_t = \sum_{s=t}^\infty \delta_s \sim \frac{\phi}{2\phi-1}\frac{1}{t}$, $t\to\infty$. From the estimate in Case 2 we know that 
\[
\begin{aligned}
0 \leq \tilde{\Delta}_t &\leq \sum_{s=T(i)}^{t-1} \delta_s \prod_{u=s+1}^{t-1} \left(1-\tfrac{1}{u}+O(\tfrac{1}{u^2})\right) 
+ \tilde{\Delta}_{T(i)} \prod_{s=T(i)}^{t-1} \left(1-\tfrac{1}{s}+O(\tfrac{1}{s^2})\right)\\
&= \Theta(1) \left[\left(\sum_{s=T(i)}^{t-1} \delta_s \tfrac{s}{t}\right) +  \tfrac{1}{t}\right],
\end{aligned}
\]
where we use that $\sum_{v=1}^t \frac{1}{v} = \log t + \gamma + O(t^{-1})$, $t\to\infty$, with $\gamma$ Euler's constant. The term between square brackets scales like $\frac{\phi}{2\phi-1}\,\frac{\log t}{t}$ as $t \to \infty$, and dominates.
\end{proof}


\subsection{Proof of Theorem~\ref{degree_growth}}
\label{ss.rate}

\begin{proof}
The proof comes in six steps: convergence, identification, upper bound, lower bound and scaling.

\medskip\noindent
{\bf 1.}
Again fix $i \geq 1$ and suppress it from the notation. By \eqref{eq:rec}, for $t \geq i$,
\[
\EE[d_{t+1} \mid \FF_t]
= d_t + \frac{\theta_t}{t(t+1)}
= d_t\left(1 + \frac{\alpha^*_t}{t+1}\right).
\]
Define 
\[
\bar M^*_t = \frac{d_t}{\prod_{s=i}^{t-1}\left(1 + \frac{\alpha^*_s}{s+1}\right)}.
\] 
Then $\EE[\bar M^*_{t+1} \mid \FF_t] = \bar M^*_t$ for all $t \geq i$, so the sequence $\{\bar M^*_t\}_{t \geq i}$ is a non-negative martingale. Hence $\lim_{t\to\infty} \bar M^*_t = \bar M^* \in [0,\infty)$ $\PP$-a.s.\ by Doob's first martingale convergence theorem. Because $\{\bar M^*>0\} \in \FF_\infty = \cup_{t \geq 1} \FF_t$, we have $\lim_{t\to\infty} \PP(\bar M^*>0 \mid \FF_t) = \1_{\{\bar M^*>0\}}$ by L\'evy's zero-one law.

\medskip\noindent
{\bf 2.}
We use the convergence established in Part 1 to identify the growth of $d_t$ as $t\to\infty$. To this end, abbreviate $\psi=1/\phi$. By Lemma \ref{lem:alphastarconv}, we have
\[
\sum_{s \geq i} \frac{|\alpha^*_s-\psi|}{\psi+s+1} < \infty \qquad \PP\text{-a.s.} 
\]
Let
\[
Q^*_t = \frac{\prod_{s=i}^{t-1} \left(1+\frac{\alpha^*_s}{s+1}\right)}{\prod_{s=i}^{t-1} 
\left(1 + \frac{\psi}{s+1}\right)}.
\] 
Since $(1+\frac{\alpha^*_s}{s+1})(1+\frac{\psi}{s+1}) =1 + \frac{\alpha^*_s-\psi}{\psi+s+1}$, it follows that $\lim_{t\to\infty} Q^*_t = Q^* \in (0,\infty)$ $\PP$-a.s. Since
\[
d_t = \bar M^*_tQ^*_t \prod_{s=i}^{t-1} \left(1 + \frac{\psi}{s+1}\right) 
\]
and $\prod_{s=i}^{t-1} (1+\tfrac{\psi}{s+1}) \sim C\,t^\psi$, $t \to \infty$, for some $C \in (0,\infty)$, it follows that $\lim_{t\to\infty} t^{-\psi} d_t = C \bar M^*Q^*$ $\PP$-a.s., which is the claim in Theorem~\ref{degree_growth} with $\epsilon_i = C \bar M^*Q^*$. Since $\bar M^*<\infty$ $\PP$-a.s., we have that $\epsilon_i <\infty$ $\PP$-a.s.

\medskip\noindent
{\bf 3.} 
To estimate $\EE[\epsilon_i]$ from above, we show that 
\begin{equation}
\label{eq:ub}
\EE[d_t] \leq c_i\,\frac{\Gamma(t)}{\Gamma(t-\psi)}, \qquad  t \geq i,
\end{equation}
with $c_i = \Gamma(i-\psi)/\Gamma(i)$, where $\Gamma$ is the Gamma-function. The proof is by induction. Clearly, the claim is true for $t=i$ because $d_i = 1$. By \eqref{eq:meanrec}, for $t \geq i$,
\[
\EE[d_{t+1}] = \EE[d_t] + \frac{\EE[\theta_t]}{t(t+1)}.
\]
Assuming that the bound in \eqref{eq:ub} holds for all $i \leq s \leq t$, we can bound
\[
\begin{aligned}
\EE[d_{t+1}] &\leq c_i\,\frac{\Gamma(t)}{\Gamma(t-\psi)} 
+ c_i\,\frac{1}{t(t+1)} \sum_{s=i}^t \frac{\Gamma(s)}{\Gamma(s-\psi)}\\
&= c_i\,\frac{\Gamma(t)}{\Gamma(t-\psi)} 
+ c_i\,\frac{1}{t(t+1)} \frac{1}{\psi+1} \left(\frac{\Gamma(t+1)}{\Gamma(t-\psi)} 
- \frac{\Gamma(i)}{\Gamma(i-1-\psi)}\right),
\end{aligned}
\]
where in the last line we use the telescoping equality
\[
\sum_{s=i}^t \frac{\Gamma(s+a)}{\Gamma(s+b)} 
= \frac{1}{a-b+1}\,\left(\frac{\Gamma(t+1+a)}{\Gamma(t+b)} - \frac{\Gamma(i+a)}{\Gamma(i-1+b)}\right), 
\qquad  a>b>-1, \quad t \geq 1,
\]
with $a=0$, $b=-\psi$. Dropping the term with the minus sign in the upper bound and putting terms together, we get
\[
\EE[d_{t+1}] \leq c_i\,\frac{\Gamma(t+1)}{\Gamma(t+1-\psi)}\,
\left(\frac{t-\psi}{t} + \frac{t-\psi}{t(t+1)}\,\frac{1}{\psi+1}\right).
\]
The term between round brackets equals $1+\chi_t$ with
\[
\chi_t = -\frac{\psi}{t}  + \frac{t-\psi}{t(t+1)}\,\frac{1}{\psi+1} = -\frac{1}{t(t+1)}, 
\]
where the last equality uses that $\psi+1=1/\psi$. Hence $\chi_t<0$ and so the bound in \eqref{eq:ub} holds for $t+1$. Thus, by induction, the bound in \eqref{eq:ub} holds for all $t \geq i$. Since $\Gamma(t)/\Gamma(t-\psi) = t^\psi[1+O(1/t)]$ as $t\to\infty$, it follows that $\limsup_{t\to\infty} t^{-\psi} \EE[d_t] \leq c_i$.

\medskip\noindent
{\bf 4.} 
To estimate $\EE[\epsilon_i]$ from below, let $\gamma_t = \EE[d_t]/(t+1)^\psi$. By \eqref{eq:meanrec}, for $t \geq i$,
\[
\gamma_{t+1} = \gamma_t\left(\frac{t+1}{t+2}\right)^{\psi}\left(1+\frac{1}{(t+1)\beta_t}\right)
\geq \gamma_t\left(\frac{t+1}{t+2}\right)^{\psi}\left(1+\frac{1}{t+1}\right)^{\frac{1}{\beta_t}}
= \gamma_t\left(\frac{t+1}{t+2}\right)^{\psi - \frac{1}{\beta_t}},
\]
where for the inequality we use that $\beta_t \geq 1$. We know from Theorem~\ref{mode_convergence} that $\lim_{t\to\infty} \beta_t = 1/\psi$, with the limit eventually reached from below. Therefore $\gamma_{t+1} > \gamma_t$ for $t$ large enough, and so $\lim_{t\to\infty} \gamma_t = \gamma > 0$ (from Part 3 we already know that $\gamma<\infty$). Iteration of the above inequality gives
\[
\frac{\gamma_t}{\gamma_i} \geq \prod_{s=i}^{t-1} \left(1-\frac{1}{s+1}\right)^{\psi-\frac{1}{\beta_s}} 
=\exp\left[-\sum_{s=i}^{t-1} \left(\psi-\frac{1}{\beta_s}\right)\,\left(\frac{1}{s} + O\left(\frac{1}{s^2}\right)\right)\right].
\]
By Lemma \ref{lem:betaconv}, the right-hand side converges to a constant $C_i > 0$ as $t\to\infty$. Moreover, $\lim_{i\to\infty} C_i =1$. Hence $\EE[\epsilon_i] \geq \gamma_i > 0$.

\medskip\noindent
{\bf 5.}
Combining Parts 3--4 with the observation that $c_i \sim i^{-\psi} \sim \gamma_i$ as $i\to\infty$, we obtain that $\lim_{t\to\infty} i^\psi \EE[\epsilon_i] = 1$. 
\end{proof}


\appendix


\section{Stochastic approximation scheme}
\label{appA}

The following definition and two lemmas are taken from \cite[Chapter 2]{B2023}.

\begin{definition}
\label{def:sas}
{\rm A recursion for a sequence $\{x_t\}_{t \geq 0}$ in $\R^d$ of the form
\begin{equation}
\label{eq:recursion}
x_{t+1} = x_t + a_t \big(h(x_t) + M_{t+1} + \epsilon_{t+1}\big), \qquad t \geq 0,
\end{equation}
is called a \emph{stochastic approximation scheme} if it satisfies the following conditions:
\begin{enumerate}
\item[(1)] 
The \emph{map} $h\colon\, \R^d \to \R^d$ is Lipschitz. 
\item[(2)] 
The \emph{step sizes} $\{a_t\}_{t \geq 0}$ satisfy $a_t \in (0,\infty)$ for all $t \geq 0$, $\sum_{t \geq 0} a_t = \infty$ and $\sum_{t\geq 0} a_t^2 < \infty$.
\item[(3)] 
The sequence $\{M_t\}_{t \geq 1}$ is a \emph{martingale difference sequence} with respect to the filtration $\{\FF_t\}_{t \geq 0}$ given by $\FF_t = \sigma (x_0, \{M_s\}_{1 \leq s \leq t})$. Moreover, the sequence $\{\zeta_t\}_{t \geq 0}$ defined by $\zeta_t = \sum_{s=0}^{t-1} a_s M_{s+1}$ is a \emph{zero-mean square-integrable martingale} that satisfies $\sum_{t \geq 0} \EE[\|\zeta_{t+1}-\zeta_t\|^2 \mid \FF_t] = \sum_{t \geq 0} a_t^2\,\EE[\|M_{t+1}\|^2 \mid \FF_t] < \infty$ $\PP$-a.s.
\item[(4)] 
The \emph{solution} $\{x_t\}_{t \geq 0}$ satisfies $\sup_{t \geq 0} \|x_t\| < \infty$ $\PP$-a.s.
\item[(5)] 
The \emph{errors} $\{\epsilon_t\}_{t \geq 1}$ are adapted to the filtration $\{\FF_t\}_{t \geq 0}$, are bounded, and satisfy $\lim_{t\to\infty} \epsilon_t = 0$ $\PP$-a.s.
\end{enumerate}
}
\end{definition}

\noindent
(Conditions (1)-(5) are common, although different conditions are possible in special cases.)  
  
The key convergence lemma for the stochastic approximation scheme is the following.
 
\begin{lemma}
\label{lem:sas}	 
Suppose that conditions {\rm (1)--(5)} in Definition \ref{def:sas} hold. Furthermore, suppose that the set of limit points of the solutions to the ODE $y'(u) = h(y(u))$, $u \geq 0$, with $y(0) \in \R^d$ coincides with the set of zeroes of $h$, and that these zeroes are isolated. Then the iterates of the recursion in \eqref{eq:recursion} converge $\PP$-a.s.\ to a (possibly sample-path dependent) zero of $h$.    
\end{lemma}

\noindent
For the case without error, i.e., $\epsilon_t = 0$ for all $t \geq 1$, Lemma~\ref{lem:sas} is derived in \cite[Section 2.1]{B2023} (see, in particular \cite[Theorem 2.1]{B2023}). The case with error can be included easily, as explained in \cite[Section 2.2]{B2023}. If $h$ has a unique zero $x^*$ that is globally attracting under the ODE, then $\lim_{t\to\infty} x_t = x^*$ $\PP$-a.s.\ for all initial values $x_0$ (see, in particular, \cite[Corollary 2.2]{B2023}).

\begin{remark}
{\rm Condition (3) is phrased differently in \cite[Chapter 2]{B2023}. Namely, there it is assumed that $\EE[\|M_{t+1}\|^2 \mid \FF_t] \leq K (1+\|x_t\|^2)$ $\PP$-a.s.\ for all $t \geq 0$ and some $K<\infty$. Together with Conditions (2) and (4) this property is used to conclude that $\sum_{t \geq 0} a_t^2\,\EE[\|M_{t+1}\|^2 \mid \FF_t] < \infty$ $\PP$-a.s., which is the key condition on which the approximation scheme is built (see \cite[Section 2.2]{B2023}). Condition (3) is the version where this key condition is included. In fact, in Section \ref{ss.fixedpointconv} we have checked Condition (3) in our particular context by using Lemma~\ref{lem:growth}.}
\end{remark}

An explicit rate of convergence estimate is derived in \cite[Chapter 2]{B2023}.  

\begin{lemma}
\label{lem:sasrate}
Put $T_0 = 0$ and $T_t = \sum_{s=0}^{t-1} a_s$, $t \geq 1$. Let $(\bar{x}(u))_{u \in [0,\infty)}$ be the piecewise linear interpolation of $\{x_t\}_{t \geq 0}$ over the successive time intervals $u \in [T_t,T_{t+1}]$, $t \geq 0$. Let $(y^{T_t}(u))_{u \in [T_t,\infty)}$ be the solution to the ODE $y'(u) = h(y(u))$, $u \geq 0$, with initial value $y^{T_t}(T_t) = \bar{x}_{T_t}$. Then, for any $t \geq 0$ and $T \in (0,\infty)$,
\begin{equation}
\label{comparison}
\sup_{u \in [T_t,T_t + T]} |\bar x(u) - y^{T_t}(u)| 
\leq K_{T,t,L}\,\eee^{LT} + C_{T,L} \sup_{t \leq s \leq s_{T,t}} a_s + \sum_{t \leq s \leq s_{T,t}} a_s |\epsilon_s|
\qquad \PP\text{-a.s.},
\end{equation}
where $s_{T,t} = \min\{s \geq t\colon\,T_s \geq T_t + T\}$ is the smallest index enveloping the time interval $[T_t,T_t + T]$, and     
\[
K_{T,t,L} = C_{T,L} L \sum_{t \leq s \leq s_{T,t}} a^2_s + \sup_{t \leq s \leq s_{T,t}} |\zeta_s-\zeta_t|, \qquad
C_{T,L} = \|h(0)\|  + L\big(C_0 + \|h(0)\| T\big)\,\eee^{LT},
\]
with $L$ the Lipschitz constant of $h$ and $C_0 = \sup_{t \geq 0} \|x_t\|$.
\end{lemma}

The estimate in \eqref{comparison} shows that the solution of the discrete stochastic approximation scheme is close to the solution of the continuous ODE for large times. For the case without error, i.e., $\epsilon_t = 0$ for all $t \geq 0$, the estimate in \eqref{comparison} is derived in \cite[Section 3.1]{B2023}. The extra term for the case with error is an obvious upper bound on the effect of the error.


\bibliographystyle{unsrt}
\bibliography{references}   


\end{document}